\input amstex
\input amsppt.sty
\magnification\magstep1

\def\ni\noindent
\def\sbs{\subset}

\def\asdim{\operatorname{asdim}}
\def\diam{\operatorname{diam}}

\def\R{\text{\bf R}}

\def\Z{\text{\bf Z}}

\def\N{\text{\bf N}}

\def\sF{\Cal F}

\def\sV{\Cal V}
\def\sU{\Cal U}
\def\sW{\Cal W}

\hoffset= 0.0in
\voffset= 0.0in
\hsize=32pc
\vsize=38pc
\baselineskip=24pt
\NoBlackBoxes
\topmatter
\author
A.N. Dranishnikov
\endauthor

\title
Groups with a polynomial dimension growth
\endtitle
\abstract
We show that finitely generated groups with a polynomial dimension growth
have Yu's property A and give an example of such groups.
\endabstract

\thanks The author was partially supported by NSF grants DMS-0305152.
\endthanks

\address University of Florida, Department of Mathematics, P.O.~Box~118105,
358 Little Hall, Gainesville, FL 32611-8105, USA
\endaddress


\email  dranish\@math.ufl.edu
\endemail

\endtopmatter

\document
\head \S1 Introduction \endhead

In the asymptotic geometry of discrete groups the growth of functions
associated to the group is of great importance. Probably the most known
concept in the area is the volume growth of a group, i.e., 
the growth of the capacity of
$r$-balls $B_r(e)$ when the radius $r$ tends to infinity. We consider
finitely generated discrete groups supplied with the word metric and look
at them as geometric objects. Following Gromov [Gr1] one can define
different geometric characteristics of a discrete group. In particular one
can speak about dimension of $r$-balls and its growth. 

First we give an informal
description of dimension.
A finite metric space $X$ of diameter $r$ can be assign a dimension on the
scale $\lambda<r$ by means of $\lambda$-approximations of $X$ by a finite
polyhedra. Under $\lambda$-approximation we assume roughly speaking
an isometric imbedding of $X$ into a regular neighborhood $N(K)$ 
(in some normed space)
of a finite polyhedron $K$ with simplexes of size $\lambda$. The minimal 
dimension of such $K$ is the dimension of $X$ on the scale $\lambda$.
If $\lambda(t)$ is a sufficiently slowly tending to infinity function, then
the growth of  the 
dimension of the $t$-ball $B_t(e)$ on the scale $\lambda(t)$ when $t$ goes 
to infinity is the dimension growth of the group $\Gamma$. We will 
give precise definitions
in the next section in the language of covers.

In contrast with the volume growth of a group all classical groups do not exhibit
any dimension growth at all. The dimension of the balls in a typical 
classical group $\Gamma$ is constant and equals
the Gromov asymptotic dimension $\asdim\Gamma$. Groups with 
finite asymptotic dimension behave nicely in the following sense:
Many famous conjecture were proved for them such as the coarse Baum-Connes and
the classical Novikov conjectures [Yu1], Gromov's hypersphericity 
conjecture [Dr2],
the $K$-theoretic integral Novikov conjecture [CG], the mod $p$ Higson conjecture
and the integral Novikov conjecture [DFW]. The asymptotic finite dimensionality
was checked for a large class of groups [Gr1],[DJ],[BD1],[BD2],[CG],[Ji],[Ro3].

In [Yu2] Guoliang Yu extended his result on the Novikov conjecture to the
groups having so called property A. In this paper we show that all finitely
generated groups with the polynomial dimension growth have property A.
Thus for groups with the polynomial dimension growth the Novikov conjecture
holds true.

Additionally, in this paper we give examples of finitely generated asymptotically 
infinite dimensional
groups with the polynomial dimension growth. 

The dimension growth like
the volume growth in the case of finitely generated groups can be at 
most exponential.
Celebrated Gromov's group [Gr2] containing an expander is an example of a group with
a nonpolynomial dimension growth. In fact one can show that the growth there
is exponential. We conclude the introduction with an open question.

\proclaim{Question}
Are there finitely generated groups of the intermediate dimension growth?
\endproclaim

\head \S2 Dimension growth of a metric space \endhead

Let $(X,\rho)$ be a metric space and let $\sU$ be a cover of $X$.
A number $\lambda>0$ is called a Lebesgue number for $\sU$ if for 
every set $A\subset X$ of diameter $\le \lambda$ there is an element
$U\in\sU$ such that $A\subset U$.
The Lebesgue number $L(\sU)$ of a cover $\sU$ is the infimum of 
all Lebesgue numbers.
The multiplicity $m(\sU)$ of a cover $\sU$ is the maximal number of
elements of $\sU$ with a nonempty intersection.
We define a function
$$
ad_{X}(\lambda)=\min\{m(\sU)\mid L(\sU)\ge\lambda\}-1.
$$
Clearly, $ad_X(\lambda)$ is monotone. We call this function the
{\it asymptotic dimension function} of a metric space $X$.

\proclaim{Proposition 2.1} For a finitely generated group $\Gamma$ with the word
metric there is $a>0$ such that $ad_{\Gamma}(\lambda)\le e^{a\lambda}$
\endproclaim
\demo{Proof}
There is $a>0$ such that $|B_{\lambda}(x)|\le e^{a\lambda}$ for all 
$x\in\Gamma$. We consider the cover $\sU$ of $\Gamma$ by the $\lambda$-balls
$\{B_{\lambda}(x)\mid x\in\Gamma\}$. Clearly, $L(\sU)\ge\lambda$. We estimate
the multiplicity $m(\sU)$ of $\sU$. If for different $x_i$ we have
$$y\in B_{\lambda}(x_1)\cap\dots\cap B_{\lambda}(x_k)$$ then
$\rho(y,x_i)\le\lambda$ for all $i$. Hence 
$k\le|B_{\lambda}(x)|\le e^{a\lambda}$. Thus,
$m(\sU)\le e^{a\lambda}$.\qed
\enddemo

The asymptotic dimension function has a direct relation to Gromov's asymptotic 
dimension.

\proclaim{Definition [Gr1]} Let $X$ be a metric space. Then 
$\asdim X\le n$ if for every $r<\infty$ there 
are uniformly bounded
$r$-disjoint families $\sU^0,\dots,\sU^n$ of subsets of $X$ such that
$\cup_i\sU^i$ is a cover of $X$.
\endproclaim

A family $\sF$ of a subsets of a metric space $(X,\rho)$ is called $r$-disjoint
if $\rho(F,F')=\inf\{\rho(x,x')\mid x\in F, x'\in F'\}\ge r$ for all 
$F,F'\in\sF$,
$F\ne F'$. We denote by $$diam(\sF)=\sup\{diam(F)\mid F\in\sF\}.$$ A family $\sF$
is called {\it uniformly bounded} if $\diam(\sF)<\infty$.

We use the notation $l_2$ for the Hilbert space of square summable sequences.
Let $\Delta$ denote the standard infinite dimensional simplex in $l_2$
$$\Delta=\{(x_1,x_2,\dots)\in l_2\mid\sum x_i=1, x_i\ge 0\}.$$
A subcomplex $K\subset\Delta$ taken with the restricted metric is called
{\it uniform}. A map $f:X\to K$ from a metric space to a uniform complex is
called {\it uniformly cobounded} if there is $D$ such that
$diam f^{-1}(St(v,K))\le D$ for all $v\in K$ where $St(v,K)$ is the star of 
a vertex $v$ in the complex $K$.
\proclaim{Proposition 2.2}
For a metric space $X$ the following conditions are equivalent:
\roster
\item{} $\asdim X\le n$;
\item{} For every $\lambda<\infty$ there is a uniformly bounded cover $\sU$ of
$X$ with $L(\sU)\ge \lambda$ and with $m(\sU)\le n+1$;
\item{} For every $\epsilon>0$ there is a uniformly cobounded $\epsilon$-Lipschitz
map $f:X\to K\subset\Delta\subset l_2$ to an $n$-dimensional uniform
polyhedron.
\endroster
\endproclaim
\demo{Proof}
(1) $\Rightarrow$ (2). Let $\lambda<\infty$ be given. Apply (1) with $r>2\lambda$
to obtain the families $\sU^0,\dots,\sU^n$. Then consider a cover
$\sU=\{N_{\lambda}(U)\mid U\in\sU^i, i=0,\dots,n\}$ where 
$N_{\lambda}(U)=\{x\in X\mid \rho(x,U)\le \lambda\}$ denotes the 
$\lambda$-neighborhood of $U$.
Since $r>2\lambda$ and each family $\sU^i$ is $r$-disjoint, we have 
$m(\sU)\le n+1$. Clearly, $L(\sU)\ge\lambda$.

(2) $\Rightarrow$ (3). Let $\sU$ be a uniformly bounded
cover of $X$ with $m(\sU)\le n+1$ and with
$L(\sU)\ge \lambda$. Then the partition of unity 
$$
\phi_U(x)=\frac{\rho(x,X\setminus U)}{\sum_{V\in\sU}\rho(x,X\setminus V)}
$$
defines the projection $p_{\sU}:X\to\Delta\subset l_2(\sU)$. It was shown in
[BD2], Proposition 1, that $p_{\sU}$ is $(2n+3)^2/L(\sU)$-Lipschitz.
The condition $m(\sU)\le n+1$ implies that $p_{\sU}(X)$ lies in an $n$-dimensional
subcomplex of $\Delta$. Note that $p_{\sU}^{-1}(St(u,K))=U$ for every vertex $u$
corresponding the element $U\in\sU$. Since $diam(\sU)<\infty$, $p_{\sU}$
is uniformly cobounded. Given $\epsilon>o$ we can take $\lambda$ large enough
that the map $p_{\sU}$ will be $\epsilon$-Lipschitz. 

(3) $\Rightarrow$ (1). Let $\sV^i=\{St(b_{\sigma},\beta^2K)\mid \sigma\subset K,
dim\sigma=i\}$, where $b_{\sigma}$ denotes the barycenter of a simplex $\sigma$ and 
$\beta^lk$ denotes the $l$-th barycentric subdivision of $K$. There are 
constants $c$ and $d$ depending only on $n$ such that the family $\sV^i$ 
is $c$-disjoint and $diam(\sV^i)\le d$ for all $i\le n$. We define
$\sU^i=f^{-1}(\sV^i)$ for a $(d/\lambda)$-Lipschitz uniformly cobounded 
map $f:X\to K$ to a uniform $n$-dimensional complex.
Then the family $\sU^i$ is $\lambda$-disjoint for all $i$. If $v\in \sigma$ is a vertex,
then $St(b_{\sigma},\beta^2K)\subset St(v,K)$. Hence 
$\diam(f^{-1}(St(b_{\sigma},\beta^2K)))\le
diam(f^{-1}(St(v,K))\le D$ and hence $\sU^i$ is uniformly bounded for every $i$.
\qed
\enddemo

\proclaim{Corollary 2.3}
The limit $\lim_{\lambda\to\infty}ad_X(\lambda)=\asdim X$.
\endproclaim

The asymptotic dimension is a {\it coarse invariant}, i.e. it is invariant 
under coarse isomorphisms. The {\it Coarse category} [Ro1-3] is a quotient
category of the category of metric spaces and coarse morphisms.
A map $f:(X,d_X)\to (Y,d_Y)$ between metric spaces is 
called {\it a coarse morphism}
if it is 
\roster
\item{} {\it metrically proper}, i.e., the preimage $f^{-1}(C)$ is bounded 
for every bounded
set $C\subset Y$;
\item{} {\it coarsely uniform}, i.e., 
$d_Y(f(x),f(x'))\le\rho(d_X(x,x'))$ for all $x,x'\in X$ for some function 
$\rho:\R_+\to\R_+$
(depending on $f$).
\endroster 
Two morphisms $f,g:X\to Y$ are equivalent if they are in
finite distance, i.e., there is
$D<\infty$ such that $d_Y(f(x),g(x))\le D$ for all $x\in X$.
The coarse category is defined as the quotient category under this equivalence 
relation. Thus, a map $f:X\to Y$ between metric spaces is a coarse isomorphism
if there is a coarse morphism $h:Y\to X$ such that $hf$ is in bounded distance
to $1_X$ and $fh$ is in bounded distance to $1_Y$.
Clearly, $\asdim X=\asdim Y$ for coarsely isometric spaces.
We note that a quasi-isometry is a coarse isomorphism.

Our main examples of metric spaces are finitely generated groups
with the left-invariant word metric.
Since for every two finite symmetric generating sets of such group
the identity map is a quasi-isometry, we obtain that the asymptotic dimension
of a finitely generated group is a well-defined group invariant.
Moreover, the growth of the asymptotic dimension function 
$ad_{\Gamma}(\lambda)$  is an invariant of a group $\Gamma$
though it is not a coarse invariant.

We recall that a morphism $f:X\to Y$ in an abstract category 
is called a {\it monomorphism} if
for every two morphisms $k,l:Z\to X$ with $f\circ k=f\circ l$ it
follows that $k=l$.
A map of metric spaces $f:X\to Y$ that represents a monomorphism in the
coarse category 
is called a {\it coarse imbedding} [Ro3]. In the literature it is often called
by an overused name a {\it uniform imbedding}.

A coarse imbedding $f:X\to Y$ defines a coarse isometry $f:X\to f(X)$.
It is not difficult to see that a map between metric spaces 
$f:X\to Y$ is a coarse embedding if and only if
there exist two monotone tending to infinity functions 
$\rho_1,\rho_2:\R_+\to\R_+$ such that
$$
(*)\ \ \ \ \ \ \ \ \ \ \rho_1(d_X(x,x'))\le d_Y(f(x),f(x'))\le \rho_2(d_X(x,x'))
$$
for all $x,x'\in X$. If $X$ is a geodesic metric space
or a group with the word metric then
the function $\rho_2$ can be chosen to be linear. 

A finitely generated
subgroup $H\subset G$ in a finitely generated group is a typical example
of a coarse imbedding. 
We say that a subgroup $H\subset G$ has a {\it polynomial distortion} [Gr1]
in $G$ if the function $\rho_1(t)$ can be chosen such that
the inverse function $\rho_1^{-1}$ is a polynomial.

For a metric space with $\asdim X\le m$ Gromov suggested to study 
the asymptotic behavior of the function [Gr1] 
$$\gamma_X(\lambda)=\inf\{diam(\sU)\mid L(\sU)\ge\lambda,\  
m(\sU)\le m+1\}$$ as the secondary dimensional invariant. 
In this paper we will refer to $\gamma_X$  as to Gromov's 
function for the inequality $\asdim X\le m$.
The asymptotic behavior of the function 
$\gamma_X(t)$ is not
a coarse invariant. It was shown in [DZ] that every
metric space $(X,d)$ with the asymptotic dimension 
$\asdim X\le n$ and with bounded
geometry is coarsely isomorphic to a metric space $(X,d')$ 
with a linear Gromov's function
for $\asdim X'\le n$. 
Nevertheless the growth of $\gamma_X$ is an invariant
of quazi-isometries and hence in the case of a finitely 
generated group $G$ the growth of 
$\gamma_G$ is an invariant of a group.

\head \S3 Property A \endhead

By $l_p$, $p\ge 1$, we denote the Banach space of sequences $\{x_n\}$ with 
the norm
$\|x\|_p=(\sum_{n=1}^{\infty}|x_n|^p)^{\frac{1}{p}}$. Also $l_p(Z)$ will denote
the corresponding Banach space with a basis indexed by a countable set 
$Z$. Thus, $l_p=l_p(\N)$.
\proclaim{Definition} Let $p\in\R_+\cup\{\infty\}$. 
A discrete metric space $(X,d)$ has property $A_p$ if there is a sequence of maps
$a^n:X\to l_p(X)$ satisfying the conditions $\|a^n_z\|_p=1$ and 
$a^n_z(y)\ge 0$ for all
$y,z\in X$, such that
\roster
\item{} there is a function $R=R(n)$ such that $supp(a^n_z)\subset B_{R}(z)$
for all $z\in X$;
\item{} for every $K>0$, $\lim_{n\to\infty}\sup_{d(z,w)<K}\|a^n_z-a^n_w\|_p=0$.
\endroster
\endproclaim
We include the case $p=\infty$ here. Recall that the norm $\|x\|_{\infty}$ is 
the $\sup$-norm.

The following proposition first was proven in [Yu2] for $p=1$.

\proclaim{Proposition 3.1} 
A discrete metric space $X$ with property $A_p$ admits a coarse 
embedding in $l_p$.
\endproclaim
\demo{Proof}
Let $a^n$ be a sequence of maps from the definition of property $A_p$.
By passing to a subsequence we may assume that 
$\sup_{d(z,w)<n}(\|a^n_z-a^n_w\|_p)^p
<1/2^n$. Let $z_0\in X$ be a base point. We define a map 
$f:X\to l_p(X\times\N)$
by the formula $f(z)(x,n)=a^n_z(x)-a^n_{z_0}(x)$. The above 
inequality insures that $f(z)\in l_p(X\times\N)$.
We shall show that $f$
is a coarse embedding.
We may assume that the function
$R$ in the definition of property $A_p$ is strictly monotone. Let $S=R^{-1}$
be the inverse function. We define $\rho_1(t)=(2S(t/2)-2)^{1/p}$ and
$\rho_2(t)=(2t+1)^{1/p}$ and check the
inequalities (*). We have
$$
(\|f(z)-f(w)\|_p)^p=\sum_{n=1}^{\infty}(\|a^n_z-a^n_w\|_p)^p\le
\sum_{n=1}^{[d(z,w)]}(\|a^n_z-a^n_w\|_p)^p+
\sum_{[d(z,w)]+1}^{\infty}(\|a^n_z-a^n_w\|_p)^p$$
$$
\le\sum_{n=1}^{[d(z,w)]}\|a^n_z-a^n_w\|_p^p+\sum_{[d(z,w)]+1}^{\infty}1/2^n
\le 2d(z,w)+1=(\rho_2(d(z,w)))^p,
$$
and
$$
(\|f(z)-f(w)\|_p)^p=\sum_{n=1}^{\infty}(\|a^n_z-a^n_w\|_p)^p\ge
\sum_{n=1}^{N}(\|a^n_z-a^n_w\|_p)^p=2N\ge(\rho_1(d(z,w)))^p.
$$
where $N=[S(d(z,w)/2)]$ is the integral part of $S(d(z,w)/2)$.
Here we used the fact that the inequality $n\le[S(d(z,w)/2)]$ implies
the inequality $R(n)\le d(z,w)/2$ which implies that $supp(a^n_z)\cap supp(a^n_w)
=\emptyset$. The latter implies that $(\|a^n_z-a^n_w\|_p)^p=(\|a^n_z\|_p)^p+
(\|a^n_w\|_p)^p=2$.\qed
\enddemo

\proclaim{Proposition 3.2}
For every finite $p\ge 1$  property $A_1$ is equivalent to  property $A_p$. 
\endproclaim
\demo{Proof} 
First we show that for $m\ge p\ge 1$  property
$A_p$ implies $A_m$.

Assume that $X$ has property $A_p$. Let $a^n:X\to l_p(X)$ be a sequence
of functions satisfying the conditions (1)-(2) from the definition of $A_p$. Then
$\sum_{y\in X}(a^n_z(y))^p=1$ and $a^n_z(y)\ge 0$ for all $z,y\in X$.
We define $b^n_z(y)=a^n_z(y)^{\frac{p}{m}}$. Clearly,  $\|b_z^n\|_m=1$. 
The condition (1)
is satisfied automatically. We check the condition (2). 
In view of the obvious inequality $t^{m/p}+(1-t)^{m/p}\le t+(1-t)=1$ 
for $m/p\ge 1$, $t\in[0,1]$,
we have that $|a-b|^{m/p}\le|a^{m/p}-b^{m/p}|$ for all $a,b\ge 0$. Hence
$$
(\|b^n_z-b^n_w\|_m)^m=\sum_{x\in X}(|b^n_z(x)-b^n_w(x)|^{m/p})^p\le
\sum_{x\in X} |b^n_z(x)^{m/p}-b^n_w(x)^{m/p}|^p=(\|a^n_z-a^n_w\|_p)^p.$$
This implies the condition (2).

Now assume that $X$ has property $A_p$, $p\in \N\setminus\{1\}$, and show
that $X$ has property $A_1$. Let $b^n:X\to l_p(X)$ be a corresponding
sequence of functions. We define $a^n=(b^n)^p$ and check that
$$
\|a^n_z-a^n_w\|_1=\sum_{x\in X} |(b^n_z(x))^p-(b^n_w(x))^p|=\sum_{x\in X} 
|b^n_z(x)-b^n_w(x)|\times|\sum_{i=0}^{p-1}(b^n_z(x))^i(b^n_w(x))^{p-1-i}|$$
By the H\"{o}lder inequality we have
$$
\le \|b_z^n(x)-b_w^n(x)\|_p\|\sum_{i=0}^{p-1}(b^n_z(x))^i(b^n_w(x))^{p-1-i}\|_q , 
$$
where $\frac{1}{p}+\frac{1}{q}=1$.
Since $pq=p+q$, we have $$(b^n_z(x))^{iq}(b^n_w(x))^{(p-1-i)q}=
(b_z^n(x)/b_w^n(x))^{iq}(b_w^n(x))^p.$$ Therefore,
$$(b^n_z(x))^{iq}(b^n_w(x))^{(p-1-i)q}\le (b_w^n(x))^p,$$ provided 
$b_z^n(x)\le b_w^n(x)$. Symmetrically we obtain that
$$(b^n_z(x))^{iq}(b^n_w(x))^{(p-1-i)q}\le (b_z^n(x))^p,$$ provided 
$b_z^n(x)\ge b_w^n(x)$. Thus,
$(b^n_z(x))^{iq}(b^n_w(x))^{(p-1-i)q}\le(b_z^n(x))^p +(b_w^n(x))^p$.
Hence $$\|(b^n_z)^i(b^n_w)^{p-1-i}\|_q=
(\sum_x(b^n_z(x))^{iq}(b^n_w(x))^{(p-1-i)q})^{1/q}\le(\sum_x(b_z^n(x))^p +
(b_w^n(x))^p)^{1/q}
=2^{1/q}.$$ Therefore,
$$
\|a^n_z-a^n_w\|_1\le 2^{1/q}p\|b_z^n-b_w^n\|_p,
$$
and the condition (2) holds for $a^n$.
\qed
\enddemo
The equality $A_1=A_2$ was proven in [Tu].
It was shown in [HR] that property $A_1$ coincides with Yu's property 
$A$ [Yu2]
for metric spaces $X$ of bounded geometry. 
Also it was proven in [HR] that
property $A$ ($=A_1=A_2$) for finitely generated groups $\Gamma$ with 
word metrics
is equivalent to the topological amenability of the natural action of $\Gamma$
on the Stone-\v{C}ech compactification $\beta\Gamma$.

\proclaim{Theorem 3.3}
Suppose that a discrete metric space $(X,\rho)$ has a polynomial 
dimension growth.
Then $X$ has property $A$.
\endproclaim
\demo{Proof}
Since $ad_X(\lambda)$ has a polynomial growth, there is $p>1$ such that
$$\lim_{\lambda\to\infty}ad_X(\lambda)/\lambda^p=0.$$ 
In view of Proposition 3.2 it suffices to prove that $X$ has property $A_p$.
Given $n<\infty$  we take an open uniformly bounded covering $\sU=\{U_i\}_{i\in J}$
of $X$ with the Lebesgue number $ L(\sU)\ge n+1$ and with the multiplicity
$m(\sU)\le ad_X(n+1)+1$. We shrink the cover $\sU$ by taking 
sets $\tilde U=N_{-n}(U)=U\setminus N_{n}(X\setminus U)$ for all $U\in\sU$
where $N_n(A)$ denotes the $n$-neighborhood of $A$. We consider
an irreducible subcover $\tilde\sV$ of this new cover, i.e., a subcover with
$$\tilde V\setminus(\bigcup_{V'\in\tilde\sV\setminus\{\tilde V\}}V')\ne\emptyset$$
for all $\tilde V\in\tilde\sV$.
Then we consider a cover $\sV=\{U\in\sU\mid \tilde U\in\tilde\sV\}$.
Note that $L(\sV)\ge n$ and $m(\sV)\le m(\sU)\le d(n+1)+1$. 
Let the cover $\sV$ be indexed by a (countable) set $J$: $\sV=\{V_j\}_{j\in J}$.
Since the cover $\tilde\sV$ is irreducible
we can define an injective map $\xi:J\to X$ such that $\xi(j)=y_j\in V_j$.
This gives us an
embedding $l_p(J)\subset l_p(X)$. Denote $\phi_i(x)=\rho(x,X\setminus U_i)$.
then for fixed $x\in X$ the family $\{\phi_i(x)\}$ defines a nonzero element
$\phi_x\in l_p(J)$. We define a map $a^n:X\to l_p(J)\subset l_p(X)$
as $a^n_z=\phi_z/\|\phi_z\|_p$ for every $z\in\Gamma$. Assume that
$\diam (\sU)\le R$. Then
$$
supp(a^n_z)=\{y_i\mid \phi_i(z)\ne 0\}=\{y_i\mid z\in U_i\in\sU\}\subset B_R(z)
$$
and the condition (1) from the definition of property $A_p$ is satisfied.
To verify the condition (2) we take $w,z\in X$ with $\rho(z,w)\le K$ and
consider the triangle inequality
$$
\|a^n_z-a^n_w\|_p\le\|a^n_z-\frac{\phi_w}{\|\phi_z\|_p}\|_p+
\|\frac{\phi_w}{\|\phi_z\|_p}-a^n_w\|_p.
$$
Since $|\rho(z,X\setminus U_j)-\rho(w,X\setminus U_j)|< \rho(z,w)$,
$\|\phi_z\|_p\ge n$, and $m(\sU)\le 2ad_X(n+1)$, we obtain the following estimate
$$
\|a^n_z-\frac{\phi_w}{\|\phi_z\|_p}\|_p=\frac{1}{\|\phi_z\|_p}\|\phi_z-\phi_w\|_p=
\frac{1}{\|\phi_z\|_p}(\Sigma_J|\rho(z,X\setminus U_j)-
\rho(w,X\setminus U_j)|^p)^{1/p}
$$
$$
\le\frac{1}{\|\phi_z\|_p}(4ad_X(n+1)\rho(z,w)^p)^{1/p}\le 4K\frac{ad_X(n+1)^{1/p}}{n}.
$$
We note that
$$
\|\frac{\phi_w}{\|\phi_z\|_p}-a^n_w\|_p=(|\frac{1}{\|\phi_z\|_p}-
\frac{1}{\|\phi_w\|_p}|)\|\phi_w\|_p=\frac{|\|\phi_z\|_p-
\|\phi_w\|_p|}{\|\phi_z\|_p}
$$
$$
\le\frac{\|\phi_z-\phi_w\|_p}{\|\phi_z\|_p}\le 4K\frac{ad_X(n+1)^{1/p}}{n}
$$
by the above estimate.
Thus $$\|a^n_z-a^n_w\|_p\le 8K\frac{ad_X(n+1)^{1/p}}{n}\to 0.$$
\qed
\enddemo

Theorem 3.3 was proved  for bounded $ad_X(\lambda)$ in [HR] and extended
for sublinear $ad_X(\lambda)$ in [Dr1].

Theorem 3.3 together with Proposition 3.2 and the result of Gouliang Yu
[Yu1] imply the following
\proclaim{Corollary 3.4}
If a group $\Gamma$ has a polynomial dimension growth then the coarse
Baum-Connes conjecture, and hence the Novikov Conjecture holds true for $\Gamma$.
\endproclaim

\proclaim{Remark 3.5}
Every discrete metric space has property $A_{\infty}$.
\endproclaim
\demo{Proof}
We define $a^n_z(x)=\max\{1-\frac{1}{n}d(x,z)),0\}$. Then
$supp(a^n_z)\subset B_{n}(z)$. Let $d(z,w)\le K$.
There is the inequality 
$$\sup_x\{1-\frac{1}{n}d(x,z)| d(x,w)>n, d(x,z)\le n\}\le K/n$$ 
since $K\ge d(z,w)\ge d(x,w)-d(x,z)\ge n-d(x,z)$.
Similarly
$$\sup_x\{1-\frac{1}{n}d(x,w))| d(x,z)>n, d(x,w)\le n\}\le K/n.$$
Note that  $\sup_x\frac{1}{n}|d(x,z))-d(x,w)|\le K/n.$
Thus we obtain that
$\|a^n_z-a^n_w\|_{\infty}\le K/n.$

Therefore the condition
(2) is satisfied.\qed
\enddemo

\head \S4 Dimension growth of wreath product \endhead

Let $G$ and $N$ be a finitely generated groups and
let $1\in G$ and $e\in N$ be their units. The {\it support} of a function 
$f:N\to G$ is the set
$$ supp(f)=\{x\in N\mid f(x)\ne 1\}.$$
The direct sum $\oplus_NG$ of groups $G$ (or restricted direct product)
is the group of functions
$C_0(N,G)=\{f:N\to G\}$ with finite supports. If $N$ is a group, there is
a natural action of $N$ on $C_0(N,G)$:  $a(f)(x)=f(xa^{-1})$ for all $a\in N$
and $f\in C_0(N,G)$. The semidirect product $C_0(N,G)\rtimes N$ is called 
a {\it restricted wreath product} and it is denoted as $N\wr G$.

REMARK. If $G$ has an element of infinite order, then $\asdim N\wr G=\infty$.
Indeed, then $N\wr G$ contains $\Z^m$ for all $m$.

We recall that the product in $N\wr G$ is defined by the formula
$$
(f,a)(g,b)=(fa(g),ab).
$$
Let $S$ and $T$ be finite generating sets for $G$ and $N$ respectively. 
Let $1\in C_0(N,G)$ 
denote
the constant function taking value $1$, and let $\delta_v^a:N\to G$ be the 
$\delta$-function, i.e., $\delta_v^a(v)=a$ and $\delta_v^a(x)=e_G$ 
for $x\ne v$. 

Then $\tilde S=\{(\delta_e^s,e),\ (1,t)\mid s\in S, t\in T\}$
is a generating set for $N\wr G$ . 
Let $\rho$ be the word metric on $N\wr G$ defined by the set $\tilde S$
and let $\|\ \|$ be the corresponding norm.

\proclaim{Proposition 4.1} Let $\pi_A:C_0(N,G)=\oplus_{x\in N}G_x\to 
\oplus_{x\in A}G_x$, 
$A\subset N$ , $G_x=G$, be
the projection. Then 
\roster
\item{} $\pi_A$ is 1-Lipschitz with respect to the metric $\rho$;
\item{} $\|z\|\ge \|a\|$ for all $z\in G_a$, $a\in N$.
\endroster
\endproclaim
\demo{Proof}
(1) We will use abbreviations $\delta_x^g$ for $(\delta_x^g,e)$ and 
$t$ for $(1,t)$ for elements of the group
$N\wr G$. Every element $w\in N\wr G$ can be presented as a word
$$
w=u_1\delta^{g_1}_eu_2\delta^{g_2}_e\dots u_n\delta^{g_n}_eu_{n+1}
$$
where $u_i\in N\setminus\{e\}$ and $g_i\in G\setminus\{1\}$
and the multiplication $u\delta$ means the action of $u$ on $\delta$. Moreover, 
every presentation of $w$
by elements of $\tilde S$ can be reduced to this one by the multiplication in 
groups $N$ and $G$.
Note that for the natural projection $p:N\wr G\to N$ we have $p(w)=
u_1\dots u_{n+1}$.
Then $C_0(N,G)$ consists of elements $w$ which have $u_1\dots u_{n+1}=e$ 
for every presentation
of the above type.
We note that the norm of $w$ being the minimal number of elements
of $\tilde S$ in a presentation of $w$ can be written as
 $$\|w\|=\min\{\sum\|u_i\|_N+\sum\|g_i\|_G\ 
 \mid w=u_1\delta^{g_1}_eu_2\delta^{g_2}_e\dots u_n\delta^{g_n}_eu_{n+1}\}.$$
Since $x\delta_e^g=\delta_x^g$, for $w\in C_0(N,G)$ we have 
$$w=u_1\delta^{g_1}_eu_2\delta^{g_2}_e
\dots u_n\delta^{g_n}_eu_{n+1}=\delta_{u_1}^{g_1}\delta_{u_1u_2}^{g_2}
\dots\delta_{u_1\dots u_n}^{g_n}.$$
Then $\pi_A(w)=\delta^{g_{i_1}}_{u_1\dots u_{i_1}}\dots
\delta^{g_{i_k}}_{u_1\dots u_{i_k}}$
where $u_1\dots u_{i_m}\in A$ for all $m\le k$. Thus, a presentation for 
$\pi_A(w)$ can be obtained by
deletion of $\delta_e^{g_r}$  for all $r$ with $u_1\dots u_r\notin A$ from
$$u_1\delta^{g_1}_eu_2\delta^{g_2}_e\dots u_n\delta^{g_n}_eu_{n+1}.$$
If one start from a shortest presentation for $w$ in the alphabet $\tilde S$, 
then after the above cancellation he
obtains a presentation for $\pi_A(w)$. It means that 
$\|w\|\ge\|\pi_A(w)\|$ for all $w\in C_0(N,G)$.

(2) For every $v\in G_a=aG_ea^{-1}$ and every shortest presentation
$$v=u_1\delta^{g_1}_eu_2\delta^{g_2}_e\dots u_n\delta^{g_n}_eu_{n+1}$$ we 
should have $u_1=a$
and hence $\|v\|\ge\|a\|$.\qed
\enddemo

\proclaim{Proposition 4.2}
Let $\pi:G\to H$ be an epimorphism of finitely generated groups
with the kernel $K=ker(\pi)$. Let $\sU$ be a cover of $H$ with
$L(\sU)\ge\lambda$ and $diam(\sU)\le R$. Let $\sV$ be a cover
of $K$ with $L(\sV)\ge 6R$ and $diam(\sV)\le D$ where we 
$K$ is taken with the 
metric restricted from $G$. Then there is
a cover $\sW$ of $G$ with $L(\sW)\ge\lambda$, $diam(\sW)\le D+2R$,
and with $m(\sW)\le m(\sU)m(\sV)$.
\endproclaim
\demo{Proof}
Let $S$ be a finite symmetric generating set for $G$ and let 
$\bar S=\pi(S)$ be a generating set for $G$. 
We denote the corresponding metrics on $G$ and $H$ by $d$ and $\rho$
respectively.
Then $\pi$ is 1-Lipschitz
for these metrics. Hence for every $R$ we have 
$\pi(N_R(K))\subset
B_R(e)$ where $N_R(A)$ denotes the $R$-neighborhood of $A$ and 
$B_R(e)$ denotes the $R$-ball centered at $e\in H$. 
Thus, $N_R(K)\subset \pi^{-1}(B_R(e))$. We shall establish the equality
$N_R(K)= \pi^{-1}(B_R(e))$.
Indeed, if $\|\pi(y)\|=k\le R$ then $\pi(y)=\pi(s_{i_1})\dots \pi(s_{i_k})$
for some $s_{i_1},\dots,s_{i_k}\in S$. Then 
$b=y(s_{i_1}\dots s_{i_k})^{-1}\in K$
and $d(y,K)\le d(y,b)=\|s_{i_1}\dots s_{i_k}\|\le k\le R$, i.e., 
$y\in N_R(K)$.

In a metric space $(X,d)$ we denote by 
$N_{-R}^X(A)=A\setminus N_R(X\setminus A)$ the $(-R)$-neighborhood of $A$.
We define a cover $\bar\sV$ of the $R$-neighborhood $N_R(K)$
as follows
$$
\bar\sV=\{N_R(N_{-2R}^{K}(V))\mid V\in\sV\}.
$$
Since $L(\sV)>5R$, we have $$L(\{N_{-2R}^{K}(V) \mid V\in\sV\})\ge 3R.$$
Then $L(\bar\sV)\ge R\ge \lambda$ where $\bar\sV$ is considered as a cover of 
$N_R(K)$. Indeed, for every $z\in N_R(K)$ we take $x\in K$ with
$d(x,z)\le R$. Then we take $V\in\sV$ with $d(x,K\setminus V)\ge 6R$.
Then $$d(x,K\setminus(N_{-2R}^K(V))\ge 4R.$$ For every $y\in B_R(z)\cap N_R(K)$
and every $y'\in K$ with $d(y,y')\le R$ we have 
$d(y,y')\le d(x,z)+d(z,y)+d(y,y')\le 3R$ and hence
$y'\in N_{-2R}^K(V)$.
Therefore $$y\in \bar V=N_R(N_{-2R}^K(V)).$$ Thus, 
$$B_R(z)\cap N_R(K)\subset\bar V.$$

Note that $$m(\bar\sV)\le m(\sV),\ \ \ \ \ diam(\bar\sV)\le diam(sV)+2R=D+2R,\ \
\text{and}\ \ L(\bar\sV)\ge R\ge \lambda.$$

For every $U\in\sU$ we fix $z_U\in G$ such that 
$\rho(\pi(z_U),H\setminus U)\ge\lambda$.
We define  
$$\sW=\{z_U\bar V\cap\pi^{-1}(U)\mid U\in\sU,\bar V\in\bar\sV\}.$$
In view of the equality $N_R(K)=\pi^{-1}(B_R(e))$ we have that $\sW$ is 
a cover of $G$.

Note that $$diam(z_U\bar V\cap\pi^{-1}(U))\le diam(z_U\bar V)=
diam(\bar V)\le D+2R.$$
We show
that $L(\sW)\ge\lambda$. Let $y\in G$ and let $U\in\sU$ be such that
$\rho((\pi(y),H\setminus U)\ge\lambda$. Since $\pi$ is 1-Lipschitz,
$$d(y,G\setminus\pi^{-1}(U))\ge\lambda.$$
Note that $$z_U^{-1}y\in z_U^{-1}\pi^{-1}(U)=\pi^{-1}(\pi(z_U)^{-1}U)
\subset\pi^{-1}(B_R(e))=N_R(K).$$
Hence there is $\bar V\in\bar\sV$ such that 
$d(z_U^{-1}y,N_R(K)\setminus\bar V)\ge\lambda$.
Hence $$d(y,z_U(N_R(K)\setminus\bar V))\ge\lambda.$$
Take $W=z_U\bar V\cap\pi^{-1}(U)\in\sW$. Then
$$G\setminus W=(z_U(N_R(K))\setminus\bar V)\ \bigcup\ (G\setminus
\pi^{-1}(U)).$$ Then $$d(y,G\setminus W)\ge
\max\{d(y,z_U(N_R(K))\setminus\bar V), d(y,G\setminus
\pi^{-1}(U))\}\ge  \lambda.$$

Clearly, $m(\sW)\le m(\sU)m(\sV)$.
\qed
\enddemo

\proclaim{Proposition 4.3}
Let $G=G_n\supset G_{n-1}\supset\dots\supset G_1\supset G_0=1$ 
be a lower central series for a finitely generated nilpotent group $G$.
Then for every $k$ the subgroup $G_{k-1}$ has a polynomial distortion in
$G_k$.
\endproclaim
\demo{Proof}
The following more general fact was cited in [Gr1] as well-known:
{\it For finitely generated nilpotent groups $H\subset G$ the subgroup $H$ 
has a polynomial distortion.}
The proof of this fact was sketched in [Gr1] for nilpotent Lie groups. 
Using the Mal'cev 
completion one can derive it
for discrete torsion free
nilpotent groups. Since every finitely generated nilpotent group is 
quasi-isometric to a finitely generated torsion free nilpotent group,
the result follows in the general case.

For the sake of completeness we give an alternative proof of the proposition.
Without loss of generality we may assume that $G$ is torsion free.
Let $(x_1,x_2,\dots,x_m)$ be
a Hall's basis for $G$ [Ha]. Thus $G_k=\langle x_1,\dots, x_{i_k}\rangle$ 
for every $k$.
Then every element $g\in G$ is uniquely expressible as 
$g=x_1^{a_1}x_2^{a_2}\dots x_n^{a_n}$, 
written symbolically as $x^{a}$, where $a_i\in \Z$, for each $i$, 
and the group operations on $G$ amount to prescribing 
polynomials $f=\colon \Z^n\times \Z^n\rightarrow \Z^n$ 
and $i=\colon \Z^n\rightarrow\Z^n$, 
where $x^{a}x^{b}= x^{f(a,b)}$ and $(x^{a})^{-1}=x^{i(a)}$. 
Thus, Hall's basis defines a bijection $\phi:G\to\Z^m$ by the rule
$\phi(g)=a$. We consider the word metric $d_k$ on $G_k$ defined by the set
$(x_1,x_2,\dots,x_{i_k})$.
Clearly, $\phi^{-1}$ is 1-Lipschitz. Since the operation laws in $G$
are polynomial there is a polynomial $p(t)$ such that the restriction
$\phi|_{B_t(e)}:B_t(e)\to \Z^m$ is $p(t)$-Lipschitz.
We denote by $\phi_k$ the restriction of $\phi$ to $G_k$: 
$\phi_k:G_k\to\Z^{i_k}$.
Let $\pi:\Z^{i_k}\to\Z^{i_{k-1}}$ be the projection. Then the map 
$\phi_k^{-1}\pi\phi_k:G_k\to G_{k-1}$ is a $G_{k-1}$-equivariant retraction.
Then $$d_{k-1}(x,x')=d_{k-1}(\phi_k^{-1}\pi\phi_k(x),
\phi_k^{-1}\pi\phi_k(x'))\le p(d_{k}(x,x'))d_{k}(x,x').$$\qed
\enddemo

\proclaim{Proposition 4.4}
Let $f:G\to H$ be an epimorphism of finitely generated groups
with a finitely generated kernel $K=ker(f)$ and let $\asdim H\le n$, 
$\asdim K\le k$. 
Suppose that the Gromov functions $\gamma_H$ and $\gamma_K$ taken for the
word metrics on $H$ and $K$ are bounded from above by a polynomial. 
Also assume that
$K$ has a polynomial distortion in $G$. 
Then Gromov's function $\gamma_G$ taken for $\asdim G\le (k+1)(n+1)$ 
has a polynomial growth.
\endproclaim
\demo{Proof}
Let $p_1(t)$ be a polynomial bound for $\gamma_H$ and
let $p_2(t)$ be a polynomial bound for $\gamma_K$.
Let $p_3(t)$ be a polynomial such that the inverse function $p_3^{-1}$ is 
defined and can 
serve as the lower bound $\rho_1$ in the inequality (*) for
the coarse imbedding $K\subset G$.
We define a polynomial $P(\lambda)=p_2(p_3(6p_1(\lambda)))+
2p_1(\lambda)$.

Let $\lambda>0$ be given.
Let $\sU$ be a cover of $H$ with $m(\sU)\le n+1$, $L(\sU)\ge\lambda$, 
and with $diam(\sU)\le p_1(\lambda)$.

Let $\sV$ be a
cover of $K$ (taken with the word metric) with $m(\sV)\le k+1$,
$L(\sV)\ge p_3(6p_1(\lambda))$, and with $diam(\sV)\le 
p_2(p_3(6p_1(\lambda)))$.
Then in the subspace metric this cover has the following properties:
$L(\sV)\ge 6p_1(\lambda)$ and $diam(\sV)\le p_2(p_3(6p_1(\lambda)))$. We apply
Proposition 4.2 to obtain a cover $\sW$ of $G$ with $L(\sW)\ge\lambda$,
$diam(W)\le p_2(p_3(6p_1(\lambda)))+2p_1(\lambda)=P(\lambda)$
and with $m(\sW)\le (n+1)(k+1)$. Thus $\gamma_G\le P(\lambda)$.
\qed
\enddemo

Clearly, every finitely generated abelian group $G$ has a linear Gromov
function $\gamma_G$.
In view of Propositions 4.3 and 4.4,
by induction we obtain the following.

\proclaim{Corollary 4.5}
For every finitely generated nilpotent group $N$  there is $k$ such that
$\asdim G\le k$ and Gromov's function
$\gamma_N$ is bounded from above by a polynomial.
\endproclaim

\proclaim{Theorem 4.6} Let $N$ be a finitely generated nilpotent group and
let $G$ be a finitely generated group with $\asdim G<\infty$. Then
the restricted wreath product $N\wr G$ has a polynomial dimension
growth.
\endproclaim
\demo{Proof}
Since $N$ is nilpotent, by Corollary 4.5 Gromov's function $\gamma_N$ has
a polynomial growth $t^k$, $k\in\N$, for the inequality
$\asdim N\le n$ for some $n$. Let $\lambda$ be given.
Then there is
an $R$-bounded cover $\sU$ of $N$ with $L(\sU)>\lambda$, $R\le \lambda^k$,  
and with
the multiplicity $m(\sU)\le n+1$. Let $\pi: N\wr G\to N$ be the
natural epimorphism. We note that $\pi$ is 1-Lipschitz with respect to
the metric $\rho$ on $N\wr G$ and the word metric on $N$ generated by $T$.

Let $r=6R$.
We consider $\oplus_{x\in B_r(e)}G_x\subset ker(\pi)=K$.
Let $\asdim G=m$. Then $\asdim G^l\le lm$ ([DJ]).
By Proposition 2.2 there is a uniformly bounded cover $\sV$ 
of $ G^{B_r(e)}=\oplus_{x\in B_r(e)}G_x$ with $m(\sV)\le|B_r(e)|m+1$ and
with $L(\sV)\ge r$. We define a cover $\tilde\sV$ of 
$$K=C_0(N,G)=\oplus_NG=\oplus_{B_r(e)}G\times
\oplus_{N\setminus B_r(e)}G$$ as follows
$$
\tilde\sV=\{V\times z\mid z\in\oplus_{N\setminus B_r(e)}G,\ V\in\sV\}.
$$
We note that for every $z\in\oplus_{N\setminus B_r(e)}G$ and every $V\in\sV$
the set $V\times z=z(V\times 1)$ is isometric to $V\times 1$. Therefore
the cover $\tilde\sV$ is uniformly bounded.
Note that $m(\tilde\sV)=m(\sV)\le |B_r(e)|(m+1)$.

We show that $L(\tilde\sV)\ge r$. Let $f\in K$, $f:N\to G$. We consider two cases.

(1) $supp(f)\subset B_r(e)$. Then there is $V\in\sV$ such that
$\rho(f,G^{B_r(e)}\setminus V)\ge r$. For every $h$ with 
$supp(h)\setminus B_r(e)\ne\emptyset$ in view of Proposition 4.1 we have 
$$\rho(f,h)=\rho(f(x),h(x))\ge\|x\|\ge r$$
for $x\in supp(h)\setminus B_r(e)\ne\emptyset$. 
Thus $$\rho (f,\oplus_{N\setminus B_r(e)}G)\ge r$$ and hence
$\rho(f,K\setminus (V\times 1))\ge r$.

(2) There is $x\in supp(f)\setminus B_r(e)$. We decompose $f=f'+\bar f$ where
$f'=f|_{B_r(e)}$ and $\bar f=f|_{N\setminus B_r(e)}$. Since $L(\sV)\ge r$,
there is $V\in\sV$ such that $$\rho(f',G^{B_r(e)}\setminus V)\ge r.$$
Then $f\in V\times \bar f\in\tilde\sV$. Note that
$$K\setminus(V\times\bar f)= ((G^{B_r(e)}\setminus V)\times\bar f)\cup
(G^{B_r(e)}\times(\oplus_{N\setminus B_r(e)}G\setminus\{\bar f\})).$$
By Proposition 4.1 $$\rho(f,(G^{B_r(e)}\setminus V)\times\bar f)\ge
\rho(f',G^{B_r(e)}\setminus V)\ge r.$$ If 
$$h\in
G^{B_r(e)}\times(\oplus_{N\setminus B_r(e)}G\setminus\{\bar f\}),$$ then
$$f|_{N\setminus B_r(e)}\ne h|_{N\setminus B_r(e)}.$$ Hence there is 
$y\in N\setminus B_r(e)$ such that $f(y)\ne h(y)$. Then 
by Proposition 4.1 $$\rho(f,h)\ge\rho(f(y),h(y))\ge\|y\|\ge r$$ for all
$$h\in G^{B_r(e)}\times(\oplus_{N\setminus B_r(e)}G\setminus\{\bar f\}).$$
Thus, $\rho(f,K\setminus(V\times\bar f))\ge r$.

We apply Proposition 4.2 to the epimorphism $\pi:N\wr G\to N$ and
to the covers $\sU$ and $\tilde\sV$ of $N$ and $K=ker(\pi)$
to obtain a uniformly bounded cover $\sW$ of $N\wr G$ with $L(\sW)\ge \lambda$
and with 
$$m(\sW)\le m(\sU)m(\tilde\sV)\le (n+1)(m+1)|B_r(e)|.$$
It is known that every finitely generated nilpotent group has a polynomial
volume growth.
Let $P(t)$ be a monotone polynomial such that $|B_r(e)|\le P(t)$. Then
$$m(\sW)\le (n+1)(m+1)P(r)=(n+1)(m+1)P(6R)\le
(n+1)(m+1)P(6\lambda^k).$$
\qed
\enddemo

\enddocument